\documentclass{article}
\usepackage[english]{babel}
\usepackage{amsmath, amsthm, amssymb}
\usepackage{setspace}
\usepackage{cite}
\usepackage{appendix}
\usepackage{listings}
\newtheorem{theorem}{Theorem}
\usepackage{abstract}
\usepackage{geometry}
\title{Gauss-Euler Primality Test}
\author{Almas Wang\thanks{Beijing, China. E-mail: alamusi18@sina.com.cn}}
\date{November 12, 2023}

\begin{document}
\maketitle

\begin{abstract}
This paper presents two efficient primality tests that quickly and accurately test all integers up to $2^{64}$.
\end{abstract}

\section{Introduction}
How to quickly and accurately test whether an integer is prime or not has always been an important problem in arithmetic. Starting with the Sieve of Eratosthenes more than 2,000 years ago, mathematicians have been trying various methods to solve this problem. In 2004, the Indian mathematician Agrawal and his two students Kayal and Saxena presented a polynomial-time algorithm~\cite{agrawal2004primes}, which is by far the most efficient deterministic algorithm. But even so, his time complexity has reached $O(\mathrm{log}^{6}\,n)$. Thus in practice the most widely used currently is a probabilistic algorithm, the Miller-Rabin primality test~\cite{rabin1980probabilistic}. Proposed by Professor Michael O. Rabin of the Hebrew University of Jerusalem, Israel. It has been more than 40 years since the algorithm was proposed, yet it remains one of the most important primality tests.

The second section of this paper presents a new primality test. Because it utilizes Gauss's law of quadratic reciprocity and Euler's criterion, it is named as Gauss-Euler primality test. Afterwards, it is compared with the Miller-Rabin primality test and it is concluded that both algorithms run with essentially the same efficiency. The third section combines the Gauss-Euler primality test with the Miller-Rabin primality test to derive a more efficient primality test, which can accurately test all integers up to $2^{64}$ with only 5+1 bases. In the final section, we will discuss primality testing for large integers.

\section{Gauss-Euler Primality Test}
This test requires the use of the following theorems~\cite{hardy1979introduction,hua2012introduction}.

\begin{theorem}
Let $p$ be an odd prime. Then
\begin{equation*}
\left ( \frac{2}{p} \right )= \left\{\begin{matrix}
1, & p\equiv \pm1 \pmod {8}\,, \\ 
-1, & p\equiv \pm3 \pmod {8}\,. 
\end{matrix}\right.
\end{equation*}
\end{theorem}

\begin{theorem}
\textbf{\upshape(Euler's criterion)} Let $p$ be an odd prime and $a$ be an integer coprime to $p$. Then
\begin{equation*}
a^{\frac{p-1}{2}}\equiv \left ( \frac{a}{p} \right )\pmod {p}\,.
\end{equation*}
\end{theorem}

\begin{theorem}
 Let $p$ be an odd prime. There are $\frac{p-1}{2}$ quadratic residues and $\frac{p-1}{2}$ quadratic non-residues in any reduced residue system.
\end{theorem}

\begin{theorem}
\textbf{\upshape(The law of quadratic reciprocity)} Let $p$ and $q$ be distinct odd prime. Then
\begin{equation*}
\left ( \frac{p}{q} \right )\left ( \frac{q}{p} \right )= \left ( -1 \right )^{\frac{p-1}{2}\frac{q-1}{2}}.
\end{equation*}
\end{theorem}

$\left ( \frac{2}{p} \right )$ in Theorem 1, $\left ( \frac{a}{p} \right )$ in Theorem 2, and $\left ( \frac{p}{q} \right )\left ( \frac{q}{p} \right )$ in Theorem 4 are Legendre symbols.

\subsection{Approach}
Using the above theorems, we can accurately test all integers below $2^{64}$ using only 1+4+2 bases. 

\subsubsection{First Step}
We will perform the Euler’s criterion test with 2 as a base. By Theorem 1 and Theorem 2, we have:\medskip\\
when $n$ is a prime of the form 8k+1 or 8k+7, then must satisfy
\begin{equation*}
2^{\frac{n-1}{2}}\equiv 1\pmod {n}\,,
\end{equation*}
when $n$ is a prime of the form 8k+3 or 8k+5, then must satisfy
\begin{equation*}
2^{\frac{n-1}{2}}\equiv -1\pmod {n}\,. 
\end{equation*}
If $n$ satisfies the above condition, we proceed to the next step. Otherwise $n$ is composite.\bigskip

For example, $341\equiv 5 \pmod {8}$. If 341 is prime, must satisfy $2^{\frac{341-1}{2}}\equiv -1 \pmod {341}$. However, we computed $2^{\frac{341-1}{2}}\not\equiv -1 \pmod {341}$. So, 341 is composite.\medskip

Another example, $561\equiv 1\pmod {8}$. If 561 is prime, must satisfy $2^{\frac{561-1}{2}}\equiv 1\pmod {561}$. We computed $2^{\frac{561-1}{2}}\equiv 1\pmod {561}$. So, we have to proceed to the next step of testing.

\subsubsection{Second Step}
Let $a= \left [ \sqrt{\dfrac{n}{i}} \right ]+j$. Where $i\in  \mathbb{Q}^{+}$, $j\in \mathbb{Z}$, the symbol $\left [ \sqrt{\dfrac{n}{i}} \right ]$ denotes the integer part of $\sqrt{\dfrac{n}{i}}$.
By Theorem 2 and Theorem 3, we have: if $n$ is an odd prime and $0<a<n$, must satisfy 
\begin{equation*}
a^{\frac{n-1}{2}}\equiv \pm 1\pmod {n}.
\end{equation*}
So, in this step, I chose to use the four numbers $\left [ \sqrt{n} \right ],$ $\left [ \sqrt{n} \right ]+1,$ $\left [ \sqrt{\frac{n}{2}} \right ],$ $\left [ \sqrt{\frac{n}{2}} \right ]+1$ as bases to perform Euler's criterion test. If $n$ passes the test, proceed to the next step. Otherwise $n$ is composite.\bigskip

For example, 
\begin{equation*}
46657\equiv 1\pmod {8}\  \mathrm{and} \   2^{\frac{46657-1}{2}}\equiv 1\pmod {46657}
\end{equation*}
\begin{equation*}
\left [ \sqrt{46657} \right ]^{\frac{46657-1}{2}}= 216^{\frac{46657-1}{2}}\equiv 1\pmod {46657}
\end{equation*}
\begin{equation*}
\left ( \left [ \sqrt{46657} \right ]+1 \right )^{\frac{46657-1}{2}}= 217^{\frac{46657-1}{2}}\equiv 1\pmod {46657}
\end{equation*}
\begin{equation*}
\left [ \sqrt{\frac{46657}{2}} \right ]^{\frac{46657-1}{2}}= 152^{\frac{46657-1}{2}}\equiv 1\pmod {46657}
\end{equation*}
\begin{equation*}
\left ( \left [ \sqrt{\frac{46657}{2}} \right ]+1 \right )^{\frac{46657-1}{2}}= 153^{\frac{46657-1}{2}}\equiv 1\pmod {46657}\,.
\end{equation*}
Therefore, 46657 passed the first and second steps of the test, and we have to proceed to the next step of testing.\medskip

Another example, $172081\equiv 1\pmod {8}$ and $2^{\frac{172081-1}{2}}\equiv 1\pmod {172081}$, passes the first step of the test. However, it failed this step of the test.
\begin{equation*}
\left [ \sqrt{172081} \right ]^{\frac{172081-1}{2}}= 414^{\frac{172081-1}{2}}\equiv 1\pmod {172081}
\end{equation*}
\begin{equation*}
\left ( \left [ \sqrt{172081} \right ]+1 \right )^{\frac{172081-1}{2}}= 415^{\frac{172081-1}{2}}\equiv 1\pmod {172081}
\end{equation*}
\begin{equation*}
\left [ \sqrt{\frac{172081}{2}} \right ]^{\frac{172081-1}{2}}= 293^{\frac{172081-1}{2}}\equiv 1\pmod {172081}
\end{equation*}
\begin{equation*}
\left ( \left [ \sqrt{\frac{172081}{2}} \right ]+1 \right )^{\frac{172081-1}{2}}= 294^{\frac{172081-1}{2}}\not\equiv \pm 1\pmod {172081}\,.
\end{equation*}
So, 172081 is composite.

\subsubsection{Third Step}
Some absolute Euler pseudoprimes (e.g., 46657) can easily pass the first and second steps of the test. In this step, we first solve this trouble. By looking at the absolute Euler pseudoprimes, we notice a feature~\cite{pinch2007absolute}:\medskip\\
if $n$ is an absolute Euler pseudoprime and $p$ is a prime and $p \nmid n$. Then   
\begin{equation*}
p^{\frac{n-1}{2}}\equiv 1\pmod {n}\,.
\end{equation*}
That is, no prime $p$ satisfies $p^{\frac{n-1}{2}}\equiv -1\pmod {n}$. We can use this feature to exclude all absolute Euler pseudoprimes.\medskip

We find the smallest prime $p$ of the form 4k+1 such that 
\begin{equation*}
\left ( \frac{n}{p} \right )= -1\,.
\end{equation*}
By Theorem 2 and Theorem 4, we have: if $n$ is prime, must satisfy 
\begin{equation*}
p^{\frac{n-1}{2}}\equiv -1\pmod {n}\,.
\end{equation*}
So, if $n$ satisfied $p^{\frac{n-1}{2}}\equiv -1\pmod {n}$, it is probable prime and must not be an absolute Euler pseudoprime. If $n$ is not satisfied, it must not be prime.\medskip

For example, $1729$ is the smallest absolute Euler pseudoprime. We find the smallest prime $p$ of the form 4k+1 such that $\left ( \frac{1729}{p} \right )= -1$. We have to compute one by one starting from the smallest prime of the form 4k+1. The first prime of the form 4k+1 is 5, the second is 13, the third is 17, and the fourth is 29...We computed $\left ( \frac{1729}{5} \right )= 1$, $\left ( \frac{1729}{13} \right )= 0$, $\left ( \frac{1729}{17} \right )= -1$. Thus, $17$ is the smallest prime $p$ of the form 4k+1 that we are looking for that satisfies the condition. By Theorem 2 and Theorem 4, if $1729$ is prime, must satisfy $17^{\frac{1729-1}{2}}\equiv -1\pmod {1729}$. However, we computed that $17^{\frac{1729-1}{2}}\not\equiv -1\pmod {1729}$. So, $1729$ is composite.\medskip

Another example, $\left ( \frac{46657}{5} \right )= -1$. If $46657$ is prime, must satisfy $5^{\frac{46657-1}{2}}\equiv -1\pmod {46657}$. We computed $5^{\frac{46657-1}{2}}\not\equiv -1\pmod {46657}$. So, $46657$ is composite. \bigskip

However, although we can test exclude absolute Euler pseudoprimes with this step, there is no guarantee that other pseudoprimes will not pass the test.\medskip

For example, $6164578258027337$ is composite (equal to $64107089\times 96160633$), but it passes all of the above tests.
\begin{equation*}
6164578258027337\equiv 1\pmod {8}\  \mathrm{and}\  2^{\frac{6164578258027337-1}{2}}\equiv 1\pmod {6164578258027337}
\end{equation*}
\begin{equation*}
\left [ \sqrt{6164578258027337} \right ]^{\frac{6164578258027337-1}{2}}\equiv -1\pmod {6164578258027337}
\end{equation*}
\begin{equation*}
\left ( \left [ \sqrt{6164578258027337} \right ]+1 \right )^{\frac{6164578258027337-1}{2}}\equiv 1\pmod {6164578258027337}
\end{equation*}
\begin{equation*}
\left [ \sqrt{\frac{6164578258027337}{2}} \right ]^{\frac{6164578258027337-1}{2}}\equiv 1\pmod {6164578258027337}
\end{equation*}
\begin{equation*}
\left ( \left [ \sqrt{\frac{6164578258027337}{2}} \right ]+1 \right )^{\frac{6164578258027337-1}{2}}\equiv -1\pmod {6164578258027337}
\end{equation*}
\begin{equation*}
\left ( \frac{6164578258027337}{5} \right )= -1\  \mathrm{and} \   5^{\frac{6164578258027337-1}{2}}\equiv -1\pmod {6164578258027337}\,.
\end{equation*}

Therefore, at this step, we have to perform the law of quadratic reciprocity test twice. That is, first find the smallest prime $p_{1}$ of the form 8k+5 (5, 13, 29, 37 ...), such that $\left ( \frac{n}{p_{1}} \right )= -1$, and then determine whether $p_{1}^{\frac{n-1}{2}}\equiv -1\left ( \mathrm{mod}\ n \right )$ is true. If it is not, $n$ is composite. If it is true, we proceed to find the smallest prime $p_{2}$ of the form 8k+1 (17, 41, 73, 89 ...),  such that $\left ( \frac{n}{p_{2}} \right )= -1$, and then determine whether $p_{2}^{\frac{n-1}{2}}\equiv -1\left ( \mathrm{mod}\ n \right )$ is true. If it is true, $n$ is prime. Otherwise $n$ is composite.\medskip

For example, $\left ( \frac{6164578258027337}{17} \right )=-1$. However, $17^{\frac{6164578258027337-1}{2}}\not\equiv -1\pmod {6164578258027337}$. So, $6164578258027337$ is composite.\vspace{1cm}

We can accurately test all integers up to $2^{64}$ by the three steps above.

\subsection{Algorithm Description}
\rule[2ex]{\textwidth}{0.4pt}
Input: an odd integer $n>1$ 
\begin{enumerate} 
    \item if ( $n\equiv 1\pmod 8$ \  or \  $n\equiv 7\pmod 8$ \  and \ $2^{\frac{n-1}{2}}\not\equiv 1\pmod n$ ) output “composite”;
    \item if ( $n\equiv 3\pmod 8$ \  or \  $n\equiv 5\pmod 8$ \  and \ $2^{\frac{n-1}{2}}\not\equiv n-1\pmod n$ ) output “composite”;
    \item if ( $\left [ \sqrt{n} \right ]^{\frac{n-1}{2}}\not\equiv 1 \pmod n$ \  or \  $\left [ \sqrt{n} \right ]^{\frac{n-1}{2}}\not\equiv n-1\pmod n$ ) output “composite”;
    \item if ( $\left ( \left [ \sqrt{n} \right ]+1 \right )^{\frac{n-1}{2}}\not\equiv 1\pmod n$ \  or \  $\left ( \left [ \sqrt{n} \right ]+1 \right )^{\frac{n-1}{2}}\not\equiv n-1\pmod n$ ) output “composite”;
    \item if ( $\left [ \sqrt{\frac{n}{2}} \right ]^{\frac{n-1}{2}}\not\equiv 1\pmod n$ \  or \  $\left [ \sqrt{\frac{n}{2}} \right ]^{\frac{n-1}{2}}\not\equiv n-1\pmod n$ ) output “composite”;
    \item if ( $\left ( \left [ \sqrt{\frac{n}{2}} \right ]+1 \right )^{\frac{n-1}{2}}\not\equiv 1\pmod n$ \  or \  $\left ( \left [ \sqrt{\frac{n}{2}} \right ]+1 \right )^{\frac{n-1}{2}}\not\equiv n-1\pmod n$ ) output “composite”;
    \item use trial division to find the smallest prime $p_{1}$ of the form 8k+5, such that $\left ( \frac{n}{p_{1}} \right )= -1$;
    \item if ( $p_{1}^{\frac{n-1}{2}}\not\equiv n-1\pmod n$ ) output “composite”;
    \item use trial division to find the smallest prime $p_{2}$ of the form 8k+1, such that $\left ( \frac{n}{p_{2}} \right )= -1$;
    \item if ( $p_{2}^{\frac{n-1}{2}}\not\equiv n-1\pmod n$ ) output “composite”;
    \item output “prime”;
\end{enumerate}
\rule{\textwidth}{0.4pt}\\
(See Appendix for the full C++ code.)

\subsection{Time Complexity Analysis}
The time complexity of the first and second steps of the Gauss-Euler primality test is $O(\mathrm{log}^{2}\,n)$. The time complexity of the third step is less easy to describe accurately. This is because I cannot determine an upper bound on the number of elementary operations performed required to find $p$ that satisfy $\left ( \frac{n}{p} \right )= -1$. But we can confirm that the size of $n$ has nothing to do with it. Moreover, by Theorem 3 we can conclude that: if $n$ coprime to $p$ and $n$ is not a square number (in the second step, we tested with base $\left [ \sqrt{n} \right ]$, excluded square numbers), the probability that $n$ is a quadratic non-residue of any prime $p$ is at least 50\%. Thus, it is usually only necessary to perform a small number of operations to find the smallest prime $p$ that satisfies $\left ( \frac{n}{p} \right )= -1$.

To summarise, the time complexity of the Gauss-Euler primality test should be $O(\mathrm{log}^{2}\,n)$. 

\subsection{Gauss-Euler Test and Miller-Rabin Test}
The Miller-Rabin primality test is currently the most widely used algorithm. Its principle is very simple and there are many related papers, so I will not repeat it here. (see Appendix for part of the C++ code).\medskip

Now let's compare the running efficiency of these two algorithms. Referring to the Gauss-Euler primality test, I use $2, 5, 17, \left [ \sqrt{n} \right ],\left [ \sqrt{n} \right ]+1,\left [ \sqrt{\frac{n}{2}} \right ],\left [ \sqrt{\frac{n}{2}} \right ]+1$ seven numbers as bases for the Miler-Rabin primality test ($5$ is the smallest prime of the form 8k+5 and $17$ is the smallest prime of the form 8k+1). Meanwhile, for the sake of comprehensiveness and fairness of the test, I chose the following four sets of numbers as test subjects.\medskip \\
Set 1: 1 - 199999 (the smallest 100,000 odd numbers),\\
Set 2: 2 - 1299709 (the smallest 100,000 primes),\\
Set 3: 1000000000000000001 - 1000000000000199999 (100,000 odd 19-digit numbers),\\
Set 4: 1000000000000000003 - 1000000000004133179 (100,000 19-digit primes).\medskip

The test results are as follows: 

\begin{table}[htbp]
\caption{Gauss-Euler Test and Miller-Rabin Test}     \centering
\begin{center}
\setlength\tabcolsep{30pt}
\renewcommand{\arraystretch}{1.4}
\begin{tabular}{c c c}
\hline
$n$      &  Miller-Rabin runtime (s) & Gauss-Euler runtime (s)   \\ \hline
Set 1    & 12.443                    & 12.462     \\  
Set 2    & 16.682                    & 16.720     \\  
Set 3    & 20.331                    & 20.234     \\  
Set 4    & 53.456                    & 53.828     \\ 
\hline
\end{tabular}
\end{center}
\end{table}
Through the above comparison, we see that both algorithms run with essentially the same efficiency.\medskip

It is worth noting that using $2, 5, 17, \left [ \sqrt{n} \right ],\left [ \sqrt{n} \right ]+1,\left [ \sqrt{\frac{n}{2}} \right ], \left [ \sqrt{\frac{n}{2}} \right ]+1$ as bases for the Miller-Rabin primality test does not guarantee 100\% accuracy up to $2^{64}$ (e.g., 33077785078626881 is a composite, equal to $105004421 \times 315013261$, but it passes the test).

\section{MR-GE Primality Test}
Next, we combine the Miller-Rabin primality test and the Gauss-Euler primality test to find a faster primality test.\medskip

The first approach I tried for this purpose was: first, perform the Miller-Rabin primality test using bases $2, \left [ \sqrt{n} \right ],\left [ \sqrt{n} \right ]+1,\left [ \sqrt{\frac{n}{2}} \right ],\left [ \sqrt{\frac{n}{2}} \right ]+1$. Second, find the smallest prime $p$ of the form 4k+1 that satisfy $\left ( \frac{n}{p} \right )=-1$, then use Euler's criterion to perform the law of quadratic reciprocity test. But there is a pseudoprime $17364052083370132981$ (equal to $548893 \times 1646677 \times 19211221$) below $2^{64}$ that would passes the test. So, I replaced some bases and tried another approach.\medskip

By Theorem 4, we have: when $p$ is a prime of the form 4k+3, and $p\neq q$, then  
\begin{equation*}
\left ( \frac{q}{p} \right )\left ( \frac{p}{q} \right )= \left\{\begin{matrix}
1 & \qquad \mathrm{if\  \mathit{q }\ is \ prime\  of\  the\  form\  4k+1}\,,\\ 
 -1& \qquad \mathrm{if\  \mathit{q }\ is \ prime\  of\  the\  form\  4k+3}\,.
\end{matrix}\right.
\end{equation*}

In the following test, we need to perform the law of quadratic reciprocity test using the prime of the form 4k+3. 

\subsection{The Test}

\rule[2ex]{\textwidth}{0.4pt}
Input: an odd integer $n>1$ 
\begin{enumerate}
    \item if ( $n=3$ or $n=5$ or $n=7$ ) output "prime";
    \item perform the Miller-Rabin primality test using five numbers $\ 2, \  \left [ \sqrt{n} \right ]-1,\  \left [ \sqrt{n} \right ]+1,\  \left [ \sqrt{\frac{n}{2}} \right ]-1,\\\  \left [ \sqrt{\frac{n}{2}} \right ]+1\ $ as bases; if $n$ fails the test, output "composite";
    \item use trial division to find the smallest prime $p$ of the form 4k+3 such that $\left ( \frac{n}{p} \right )= -1$;
    \item if ( $n\equiv 1\pmod 4$ \  and \ $p^{\frac{n-1}{2}}\not\equiv n-1\pmod n$ ) output “composite”;
    \item if ( $n\equiv 3\pmod 4$ \  and \ $p^{\frac{n-1}{2}}\not\equiv 1\pmod n$ ) output “composite”;
    \item output “prime”;
\end{enumerate}
\rule{\textwidth}{0.4pt}\\
(See Appendix for part of the C++ code.)\medskip

This way, we can accurately test all integers up to $2^{64}$ with just 5+1 bases.

\subsection{MR-GE Test and Gauss-Euler Test}
Now, let's compare the running time of the Gauss-Euler primality test and the MR-GE primality test. Select the following four sets of numbers again for testing.\medskip \\
Set 1: 1 - 199999 (the smallest 100,000 odd numbers),\\
Set 2: 2 - 1299709 (the smallest 100,000 primes),\\
Set 3: 1000000000000000001 - 1000000000000199999 (100,000 odd 19-digit numbers),\\
Set 4: 1000000000000000003 - 1000000000004133179 (100,000 19-digit primes).\medskip

The test results are as follows: 

\begin{table}[htbp]
\caption{MR-GE Test and Gauss-Euler Test}     \centering
\begin{center}
\setlength\tabcolsep{30pt}
\renewcommand{\arraystretch}{1.4}
\begin{tabular}{c c c}
\hline
$n$      & Gauss-Euler runtime (s)   & MR-GE runtime (s)       \\ \hline
Set 1    & 12.462                    & 12.087       \\  
Set 2    & 16.720                    & 15.653       \\  
Set 3    & 20.234                    & 19.385       \\  
Set 4    & 53.828                    & 48.402       \\ 
\hline
\end{tabular}
\end{center}
\end{table}
By comparison, we see that the MR-GE primality test takes less time to run due to the reduction of one base, and the larger $n$ is, the more obvious the effect is.

\section{Discussion}
Currently, there is no efficient deterministic algorithm that can quickly perform primality tests on large integers. However, large primes play a non-negligible role in many modern cryptosystems (e.g., RSA, ECC)~\cite{mahto2018performance}. Therefore, it is necessary for us to discuss the primality test for large integers. 

\subsection{256-Bit Primes}
The use of 256-bit primes is required in some encryption algorithms (e.g., ECC). Therefore, I hope we can find a deterministic primality test up to $2^{256}$. This is possible if we add a certain number of bases to the MR-GE primality test. To this end, we may need to spend a lot of energy and time trying various approaches. I have limited time and energy alone, and I hope that like-minded people will join me in achieving this goal.\\
Here is a approach I recommend(10+2 bases):\\
\rule[-5ex]{\textwidth}{0.4pt}
\begin{enumerate}
\item \begin{spacing}{1.25}Perform the Miller-Rabin primality test using ten numbers $\ 2,\  \left [ \sqrt{n} \right ]-1,\  \left [ \sqrt{n} \right ],\  \left [ \sqrt{n} \right ]+1,\  \left [ \sqrt{\frac{n}{2}} \right ]-1,\\\  \left [ \sqrt{\frac{n}{2}} \right ],\   \left [ \sqrt{\frac{n}{2}} \right ]+1,\  \left [ \sqrt{\frac{n}{3}} \right ]-1,\  \left [ \sqrt{\frac{n}{3}} \right ],\  \left [ \sqrt{\frac{n}{3}} \right ]+1\ $ as bases. 
\item \end{spacing}Find, in turn, the smallest prime $p_{1}$ of the form 4k+1, and the smallest prime $p_{2}$ of the form 4k+3 that satisfy $\left ( \frac{n}{p} \right )=-1$. Then use Euler's criterion to perform the law of quadratic reciprocity test.
\end{enumerate}
\rule[5ex]{\textwidth}{0.4pt}\medskip

If we find that any of the composite numbers pass the test, we need to replace or add the base and try again. In short, we still have to do a lot of work to achieve this goal.\vspace{1cm}

\subsection{2048-Bit Large Primes}
In the RSA cryptosystem, primality testing of 1024-bit or even 2048-bit large integers is required. For this purpose, we may try the following approach: (26+4 bases):\\
\rule[-5ex]{\textwidth}{0.4pt}
\begin{enumerate} 
\item 
\begin{spacing}{1.5}
Perform the Miller-Rabin primality test using  twenty-six numbers $\ 2,\  \left [ \sqrt{n} \right ]-2,\  \left [ \sqrt{n} \right ]-1,\  \left [ \sqrt{n} \right ],\  \left [ \sqrt{n} \right ]+1,\\ \left [ \sqrt{n} \right ]+2,\  \left [ \sqrt{\frac{n}{2}} \right ]-2,\  \left [ \sqrt{\frac{n}{2}} \right ]-1,\  \left [ \sqrt{\frac{n}{2}} \right ],\  \left [ \sqrt{\frac{n}{2}} \right ]+1,\  \left [ \sqrt{\frac{n}{2}} \right ]+2,\  \left [ \sqrt{\frac{n}{3}} \right ]-2,\  \left [ \sqrt{\frac{n}{3}} \right ]-1,\  \left [ \sqrt{\frac{n}{3}} \right ],\  \left [ \sqrt{\frac{n}{3}} \right ]+1,\\ \left [ \sqrt{\frac{n}{3}} \right ]+2,\  \left [ \sqrt{\frac{n}{5}} \right ]-2,\  \left [ \sqrt{\frac{n}{5}} \right ]-1,\  \left [ \sqrt{\frac{n}{5}} \right ],\  \left [ \sqrt{\frac{n}{5}} \right ]+1,\  \left [ \sqrt{\frac{n}{5}} \right ]+2,\  \left [ \sqrt{\frac{n}{7}} \right ]-2,\  \left [ \sqrt{\frac{n}{7}} \right ]-1,\  \left [ \sqrt{\frac{n}{7}} \right ],\  \left [ \sqrt{\frac{n}{7}} \right ]+1,\\ \left [ \sqrt{\frac{n}{7}} \right ]+2\ $ as bases.
\item\end{spacing} Find, in turn, the smallest prime $p_{1}$ of the form 8k+1, the smallest prime $p_{2}$ of the form 8k+3, the smallest prime $p_{3}$ of the form 8k+5, and the smallest prime $p_{4}$ of the form 8k+7 that satisfy $\left ( \frac{n}{p} \right )=-1$. Then use Euler's criterion to perform the law of quadratic reciprocity test.

\end{enumerate}
\rule[5ex]{\textwidth}{0.4pt}\medskip

I cannot guarantee 100\% accuracy of the above approach. However, if anyone discovers a composite number that passes the test, I am willing to offer a reward of \$1,000 for your efforts.

\bibliographystyle{abbrv}  
\bibliography{Gauss-Euler}

\newpage

\appendix
\section*{Appendix}
\subsection*{1. Full C++ Code for the Gauss-Euler Primality Test}
\lstset{
numbers=none,
language={[Visual]C++},
frame=shadowbox,
basicstyle=\scriptsize
}
\begin{lstlisting}
#include <iostream>
#include <cstdio>
#include <cmath> 
using namespace std;
typedef unsigned long long ll;
ll Mod_Mul(ll a, ll b, ll c) {
	ll res = 0;
	while (b) {
		if (b & 1)
		    res = (res + a) % c;
		a = (a + a) % c;
		b >>= 1;
	}  return res;
}
ll Power_Mod(ll a, ll b, ll c) {
	ll res = 1;
	while (b) {
		if (b & 1)
		    res = Mod_Mul(res, a, c);
		a = Mod_Mul(a, a, c);
		b >>= 1;
	}  return res;
}
bool Gauss_Euler(ll n) {
	if (n == 2)  return true;
	if (!(n & 1) || n < 2)  return false;
	ll t = Power_Mod(2, (n - 1) / 2, n);
	if ((n % 8 == 1 || n % 8 == 7) && t != 1) { return false;
        }
	if ((n % 8 == 3 || n % 8 == 5) && t != n - 1) { return false;
        }
	for (int j = 1; j <= 2; j++) {
		ll a = (ll)sqrt(n / j);
		for (ll i = a; i <= a + 1; i++) {
			ll q = Power_Mod(i, (n - 1) / 2, n);
			if (q != 1 && q != n - 1)  return false;
		}
	}
	ll p1, p2, i, j;
	for (p1 = 5; ; p1 += 8) {
		for (i = 3; i * i < p1; i += 2) {
			if (p1 % i == 0)  break;
		}
		for (j = 1; j <= (p1 - 1) / 2; j++) {
			if (i * i < p1 || n % p1 == 0 || n % p1 == j * j % p1)  break;
		}
		if (i * i > p1 && j > (p1 - 1) / 2)  break;
	}
	if (Power_Mod(p1, (n - 1) / 2, n) != n - 1) {
		return false;
	}
	for (p2 = 17; ; p2 += 8) {
		for (i = 3; i * i < p2; i += 2) {
			if (p2 % i == 0)  break;
		}
		for (j = 1; j <= (p2 - 1) / 2; j++) {
			if (i * i < p2 || n % p2 == 0 || n % p2 == j * j % p2)  break;
		}
		if (i * i > p2 && j > (p2 - 1) / 2)  break;
	}
	if (Power_Mod(p2, (n - 1) / 2, n) != n - 1) {
		return false;
	}
	return true;
}
int main() {
	ll n;
	while (scanf("%llu", &n) != EOF) {
		if (Gauss_Euler(n)) printf("Yes\n");
		else printf("No\n");
	}
	return 0;
}
\end{lstlisting}
(If $n$ is greater than $2^{63}$ an overflow occurs, requiring the use of \_\_int128.)

\subsection*{2. Part of the C++ Code for the Miller-Rabin Primality Test}
\begin{lstlisting}
bool Miller_Rabin(ll n){
    ll i, j, k;
    ll s = 0, t = n - 1;
    if (n == 2 || n == 5 || n == 17)  return true;
    if (n < 2 || !(n & 1))  return false;
    while (!(t & 1)) {  
        s++;
        t >>= 1;
    }
    ll m = (ll)sqrt(n);
    ll r = (ll)sqrt(n / 2);
    ll prime[7] = { 2, 5, 17, m, m + 1, r, r + 1 };
    for (i = 0; i < 7; i++){
        ll a = prime[i];
        ll b = Power_Mod(a, t, n);      
        for (j = 1; j <= s; j++) {    
            k = Mod_Mul(b, b, n);   
            if (k == 1 && b != 1 && b != n - 1)  return false;
            b = k;
        }
        if (b != 1)  return false;  
    }
    return true;
}
\end{lstlisting}\vspace{1cm}

\subsection*{3. Part of the C++ Code for the MR-GE Primality Test}
\begin{lstlisting}
bool MR_GE(ll n) {
    ll x, y, k;
    ll s = 0, t = n - 1;
    if (n == 2 || n == 3 || n == 5 || n == 7)  return true;
    if (n < 2 || !(n & 1))  return false; 
    while (!(t & 1)) {
        s++;
        t >>= 1;
    }
    ll m = (ll)sqrt(n);
    ll r = (ll)sqrt(n / 2);
    ll prime[5] = { 2, m + 1, m - 1, r + 1, r - 1 };
    for (x = 0; x < 5; x++) {
        ll a = prime[x];
        ll b = Power_Mod(a, t, n);     
        for (y = 1; y <= s; y++) {
            k = Mod_Mul(b, b, n);   
            if (k == 1 && b != 1 && b != n - 1)  return false;
            b = k;
        }
        if (b != 1)  return false;  
    }
    ll p, i, j, d;
    for (p = 3; ; p += 4) {
        for (i = 3; i * i < p; i += 2) {
            if (p % i == 0) break;
        }
        for (j = 1; j <= (p - 1) / 2; j++) {
            if (i * i < p || n % p == 0 || n % p == j * j % p) break;
        }
        if (i * i > p && j > (p - 1) / 2) break;
    } 
    d = Power_Mod(p, (n - 1) / 2, n);
    if ((n % 4 == 1 && d != n - 1) || (n % 4 == 3 && d != 1)) {
        return false;
    }
    return true;
}
\end{lstlisting}

\end{document}